\documentclass[12pt]{amsart}%{article}
\usepackage[foot]{amsaddr}
\usepackage[utf8]{inputenc}
\usepackage{mathrsfs, microtype}
\usepackage{fullpage}
\usepackage{amssymb}
\usepackage{amsmath}
\usepackage{epsfig}
\usepackage{algorithm}
\usepackage{algorithmic}
\usepackage{latexsym}
\usepackage{amsmath}
\usepackage{amsfonts}
\usepackage{graphicx}
\usepackage[normalem]{ulem}
\usepackage{array}
\usepackage{colortbl}
\usepackage{bm}
\usepackage{color}\usepackage[dvipsnames]{xcolor}
\usepackage{tikz}
\usetikzlibrary{patterns, calc}
\usepackage{array,multirow}%,makecell}
\usepackage[
  bookmarks=false, % with bookmarks \protect\tikz doesn't  work inside \section{}
  colorlinks,
  citecolor=brown!70!black,
  linkcolor=brown!80!black,
  urlcolor=blue!70!black,
]{hyperref}
\usepackage{relsize,exscale}
    \newenvironment{dedication}
        {\vspace{-1ex}\begin{quotation}\begin{center}\begin{em}}
        {\par\end{em}\end{center}\end{quotation}}

\newtheorem{lem}{Lemma}
\newtheorem{defn}{Definition}
\newtheorem{cor}{Corollary}
\newtheorem{thm}{Theorem}
\newtheorem{rem}{Remark}

\newtheorem{prop}{Proposition}
\newtheorem{conj}{Question}

\newcommand{\comm}[1]{}
 	\definecolor{lightlightgray}{rgb}{0.93, 0.93, 0.93}
 		\definecolor{llightgray}{rgb}{0.87, 0.87, 0.87}

\newcolumntype{C}[1]{>{\centering\arraybackslash }b{#1}}

\newcommand{\s}{{\texttt s}}
\newcommand{\st}{{\texttt t}}

\newcommand\oeis[1]{\href{https://oeis.org/#1}{#1}}

\begin{document}

\title{A lattice on Dyck paths  close to the Tamari lattice}

\author{Jean-Luc Baril$^1$}
\author{Sergey Kirgizov$^1$}
\address{$^1$\rm LIB, Universit\'e de Bourgogne,
  B.P. 47 870, 21078 Dijon Cedex France}
\email{barjl@u-bourgogne.fr}
\email{sergey.kirgizov@u-bourgogne.fr}

\author{Mehdi Naima$^2$}
\address{$^2$\rm Sorbonne Université, CNRS, LIP6, F-75005 Paris, France}
\email{mehdi.naima@lip6.fr}
\date{}

\maketitle

\begin{dedication}
In memory of Jean-Marcel Pallo.
\end{dedication}

\vspace{1em}
\noindent\textbf{Abstract.}
\quad % Indent first line of paragraph
We introduce a new poset structure on Dyck paths where the covering
relation is a particular case of the relation inducing the Tamari
lattice. We prove that the transitive closure of this relation
endows Dyck paths with a lattice structure. We provide a trivariate
generating function counting the number of Dyck paths with respect to
the semilength, the numbers of outgoing and incoming edges in the
Hasse diagram. We deduce the numbers of coverings, meet and join
irreducible elements. As a byproduct, we present a new involution
on Dyck paths that transports the bistatistic of the numbers of outgoing and incoming edges to its reverse. Finally, we give a generating function for the number of
intervals, and we compare this number with the number of intervals in the Tamari lattice.

\vspace{0.5cm}
\noindent\textbf{Keywords:}
% \quad % Indent first line of paragraph
Enumeration, Dyck path, Tamari lattice, join and meet irreducible elements, interval enumeration, involution.
%\subjclass{05A15, 05A19.}

\section{Introduction and motivation}

Many various classes of combinatorial objects are enumerated by the
well-known Catalan numbers. For instance, it is the case of Dyck
paths, planar trees, rooted binary trees, triangulations, Young
tableaux, non-associative products, stack sortable permutations,
permutations avoiding a pattern of length three, and so on. A list of
over 60 types of such combinatorial classes of independent interest
has been compiled by Stanley~\cite{Sta}. Generally, these classes have
been studied in the context of the enumeration according to the length
and given values of some parameters. Many other works investigate structural
properties of these sets from order theoretical point of view. Indeed,
there exist several partial order relations on Catalan sets which
endow them  with an interesting lattice structure
\cite{Barc,BMKN,Barp1,Barp2,Pal1,Simi,Tam}. Of much interest is probably
the so-called Tamari lattice~\cite{Huan,Tam} which can be obtained equivalently in
different ways. The coverings of the Tamari lattice could be different
kinds of elementary transformations as reparenthesizations of letter
products~\cite{Grat}, rotations on binary trees~\cite{Pal2,Slea}, diagonal
flips in triangulations~\cite{Slea}, and rotations on Dyck paths
\cite{Barmo,Berg,Chapo,Pal2}. The Tamari lattice appears in different domains. Its Hasse
diagram is a graph of the polytope called associahedron; it is a
Cambrian lattice underlying the combinatorial structure of Coxeter
groups; and it has many enumerative properties with tight links with
combinatorial objects such as planar maps~\cite{Bous, Boucha,Chap,Fan,Mul}.

In this paper, we introduce a new poset structure on Dyck paths where
the covering relation is a particular case of the covering relation that generates the Tamari lattice. This is a new attempt to
study a lattice structure on a Catalan set, and since it is close to the Tamari lattice, its study becomes natural.  We prove that the
transitive closure of this relation endows Dyck paths with a lattice structure. It is worth noting that this lattice has been referred to as the {\it pyramid lattice} in \cite{BMKN}. We provide a trivariate generating function for the number of Dyck paths with respect to the semilength, and the statistics $\s$ and $\st$ giving the number of outgoing edges and the number of incoming edges, respectively. As a byproduct, we obtain generating functions and closed forms for the numbers of meet and join irreducible elements, and for the number of coverings. An asymptotic is given for the ratio between the numbers of coverings in the Tamari lattice and those in our lattice. Also, we exhibit an involution on the set of Dyck paths that transports the bistatistic $(\s,\st)$ to $(\st,\s)$. Finally, we provide the generating function for the number of intervals and we offer some open problems.

\section{Notation and definitions}

In this section, we provide necessary notation and  definitions in the context of Dyck paths, combinatorics and order theory.

\begin{defn}
  A Dyck path is a lattice path in $\mathbb{N}^2$ starting at the
  origin, ending on the $x$-axis and consisting of up-steps $U=(1,1)$
  and down-steps $D=(1,-1)$.
\end{defn}

Let $\mathcal{D}_n$ be the set of Dyck paths of semilength $n$ (i.e.,
with $2n$ steps), and $\mathcal{D}=\bigcup_{n\geq 0}\mathcal{D}_n$. The
cardinality of $\mathcal{D}_n$ is the $n$-th Catalan number $c_n =
(2n)!/(n!(n+1)!)$ (see~\oeis{A000108} in~\cite{Sloa}). For instance,
the set $\mathcal{D}_3$ consists of the five paths $UDUDUD$, $UUDDUD$,
$UDUUDD$, $UUDUDD$ and $UUUDDD$.

In this work, we will use the {\it first return decomposition} of a
Dyck path, $P=U R D S$, where $R,S\in\mathcal{D}$, and the {\it last
  return decomposition} of a Dyck path, $P=R U SD$, where
$R,S\in\mathcal{D}$. Such a decomposition is unique and will be used
to obtain a recursive description of the set $\mathcal{D}$. A Dyck
path having a first return decomposition with $S$ empty will be called
${\it prime}$, which means that the path touches the $x$-axis only at
the origin and at the end.

\begin{figure}[h]
\begin{center}($i$)\quad\scalebox{0.55}{\begin{tikzpicture}[ultra thick]
 \draw[black, thick] (0,0)--(8,0); \draw[black, thick] (0,0)--(0,2.5);
  \draw[black, line width=1.5pt] (0,0)--(0.8,0.8)--(1.2,0.4)--(2,1.2)--(2.4,0.8);
  \draw[orange, line width=3pt] (2.4,0.8)--(2.8,0.4)--(3.6,1.2)--(4,0.8)--(4.4,1.2)--(5.2,0.4);
  \draw[black, line width=1.5pt] (5.2,0.4)--(5.6,0.8)--(6.4,0);
 \end{tikzpicture}}\quad $\stackrel{\mathcal{T}}{\longrightarrow}$\quad\scalebox{0.55}{\begin{tikzpicture}[ultra thick]
  \draw[black, thick] (0,0)--(8,0); \draw[black, thick] (0,0)--(0,2.5);
  \draw[black, line width=1.5pt] (0,0)--(0.8,0.8)--(1.2,0.4)--(2,1.2)--(2.4,0.8);
  \draw[orange, line width=3pt] (2.4,0.8)--(3.2,1.6)--(3.6,1.2)--(4,1.6)--(5.2,0.4);
 \draw[black, line width=1.5pt] (5.2,0.4)--(5.6,0.8)--(6.4,0);
 \end{tikzpicture}}
\end{center}
\begin{center}($ii$)\quad \scalebox{0.55}{\begin{tikzpicture}[ultra thick]
 \draw[black, thick] (0,0)--(8,0); \draw[black, thick] (0,0)--(0,2.5);
  \draw[black, line width=1.5pt] (0,0)--(0.8,0.8)--(1.2,0.4)--(2,1.2)--(2.4,0.8);
  \draw[orange, line width=3pt] (2.4,0.8)--(2.8,0.4)--(3.6,1.2)--(4,1.6)--(4.4,1.2)--(5.2,0.4);
  \draw[black, line width=1.5pt] (5.2,0.4)--(5.6,0.8)--(6.4,0);
 \end{tikzpicture}}\quad $\longrightarrow$\quad\scalebox{0.55}{\begin{tikzpicture}[ultra thick]
  \draw[black, thick] (0,0)--(8,0); \draw[black, thick] (0,0)--(0,2.5);
  \draw[black, line width=1.5pt] (0,0)--(0.8,0.8)--(1.2,0.4)--(2,1.2)--(2.4,0.8);
  \draw[orange, line width=3pt] (2.4,0.8)--(3.2,1.6)--(3.6,2)--(4,1.6)--(5.2,0.4);
 \draw[black, line width=1.5pt] (5.2,0.4)--(5.6,0.8)--(6.4,0);
 \end{tikzpicture}}
\end{center}
\caption{ ($i$) corresponds to a covering relation for the Tamari Lattice, while ($ii$) corresponds to the covering relation for the new lattice of this study.}
\label{fig1}\end{figure}
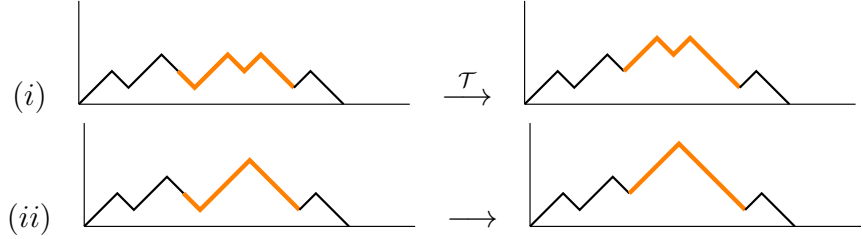

\begin{defn}
  A {\it peak} in a Dyck path is an occurrence of the subpath $UD$. A
  {\it pyramid} is a maximal occurrence of $U^kD^\ell$, $k,\ell\geq
  1$, in the sense that this occurrence cannot be extended in a
  occurrence of $U^{k+1}D^\ell$ or in a occurrence of
  $U^{k}D^{\ell+1}$.
\end{defn}
We say that a pyramid $U^kD^\ell$ is {\it symmetric} whenever
$k=\ell$, and {\it asymmetric} otherwise. The {\it weight} of a
symmetric pyramid $U^kD^k$ is $k$. For instance, the Dyck path in the
south west of Figure~\ref{fig1} contains two symmetric pyramids and two
asymmetric pyramids.

\medskip

A {\it statistic} on the set $\mathcal{D}$ of Dyck paths is a function
$\texttt{f}$ from $\mathcal{D}$ to $\mathbb{N}$, and a multistatistic
is a tuple of statistics $(\texttt{f}_1,\texttt{f}_2, \ldots,
\texttt{f}_t)$, $t\geq 2$. Given two statistics (or multistatistics)
$\texttt{f}$ and $\texttt{g}$, we say that they have the same {\it
  distribution} (or equivalently, are {\it equidistributed}) if, for
any $k\geq 0$, $$\texttt{card}\{P\in\mathcal{D},
\texttt{f}(P)=k\}=\texttt{card}\{P\in\mathcal{D}, \texttt{g}(P)=k\}.$$

Below, we define two important statistics for our study.

\begin{defn}
  Let $\s$ be the statistic on $\mathcal{D}$ where $\s(P)$ is the
  number of occurrences in $P$ of $DU^kD^k$, $k\geq 1$. Let $\st$ be
  the statistic on $\mathcal{D}$ where $\st(P)$ is the number of
  occurrences in $P$ of $U^kD^kD$, $k\geq 1$.
\label{def3}
\end{defn}
 For instance, for the path $P$ in the south west of
 Figure~\ref{fig1}, we have $\s(P)=3$ and $\st(P)=4$.
 \medskip

We end this section by defining the main concepts of order theory that
we will use in this paper. We can find all these definitions in
\cite{Grat} for instance.

A {\it poset} $\mathcal{L}$ is a set endowed with a partial order
relation. Given two elements $P,Q\in \mathcal{L}$, a {\it meet} (or
{\it greatest lower bound}) of $P$ and $Q$, denoted $P\wedge Q$, is an
element $R$ such that $R\leq P$, $R\leq Q$, and for any $S$ such that
$S\leq P$, $S\leq Q$ then we have $S\leq R$. Dually, a {\it join} (or
{\it least upper bound}) of $P$ and $Q$, denoted $P\vee Q$, is an
element $R$ such that $P\leq R$, $Q\leq R$, and for any $S$ such that
$P\leq S$, $Q\leq S$ then we have $R\leq S$. Notice that join and meet
elements do not necessarily exist in a poset. A {\it lattice} is a
poset where any pair of elements admits a meet and a join.

\begin{defn}
An element $P\in\mathcal{L}$ is {\it join-irreducible} (resp., {\it
  meet-irreducible}) if $P=R\vee S$ (resp., $P=R\wedge S$) implies
$P=R$ or $P=S$.
\end{defn}

\begin{defn}
  An {\it interval} $I$ in a poset $\mathcal{L}$ is a subset of
  $\mathcal{L}$ such that for any $P,Q\in I$, and any $R\in
  \mathcal{L}$, if $P\leq R$ and $R\leq Q$, then $R$ is also in $I$.
\end{defn}

In 1962~\cite{Tam}, the Tamari lattice $\mathcal{T}_n$ of order $n$ is
defined by endowing the set $\mathcal{D}_n$ with the transitive
closure $\preceq$ of the covering relation $$P\stackrel{\mathcal{T}}{\longrightarrow} P'$$
that transforms an occurrence of $DUQD$ in $P$ into an occurrence $UQDD$ in $P'$, where
$Q$ is a Dyck path (possibly empty). The top part of Figure~\ref{fig1} shows
an example of such a covering, and Figure~\ref{fig2} illustrates the
Hasse diagram of $\mathcal{T}_n$ for $n=4$ (the red edge must be
considered). The number of meet (resp., join) irreducible elements is
$n(n-1)/2$, and the number of coverings is $(n-1) c_n/2$~\cite{Gey}
where $c_n$ is the $n$-th Catalan number. In 2006,
Chapoton~\cite{Chap} proved that the number of intervals in
$\mathcal{T}_n$ is $$\frac{2}{n(n+1)}\binom{4n+1}{n-1}.$$
Later, Bernardi and Bonichon~\cite{Bern} exhibited a bijection between intervals in $\mathcal{T}_n$ and  minimal realizers.

Now, we introduce a new partial order on $\mathcal{D}_n$ for $n\geq 0$.  We endow it with the order relation $\leq$ defined by the transitive closure of the  covering relation  $$P\longrightarrow P'$$
that transforms an occurrence of $DU^kD^k$ in $P$ into an occurrence $U^kD^kD$ in $P'$, where $k\geq 1$.
For short, we will often use the notation
$$DU^kD^k\longrightarrow  U^kD^kD, ~k\geq 1$$
whenever we need to show where the transformation is applied.

\begin{figure}[ht]
     \centering
     \includegraphics[scale = 0.5]{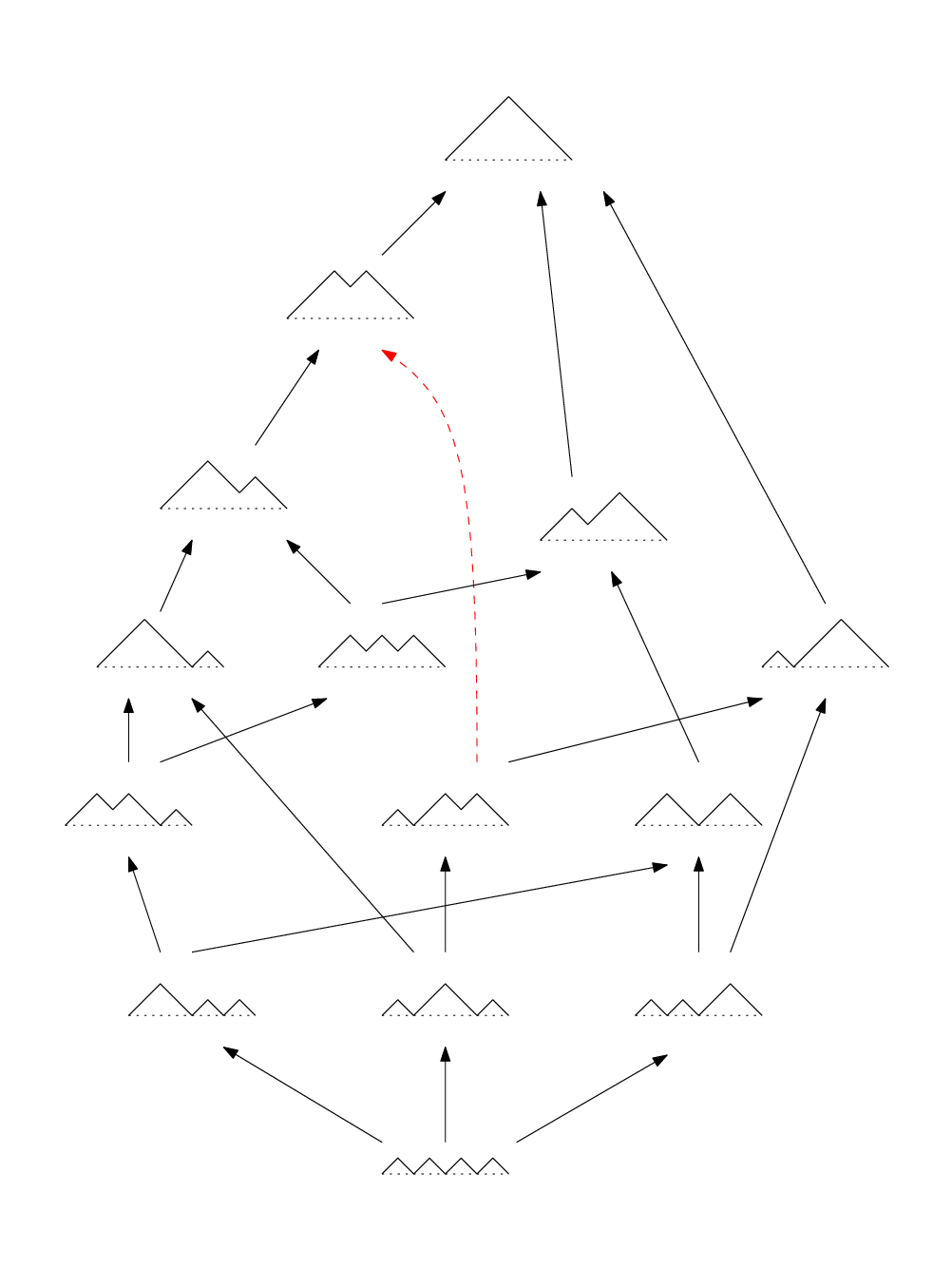}
     \caption{The Hasse diagram of $\mathcal{S}_4=(\mathcal{D}_4,\leq)$. The Tamari lattice $\mathcal{T}_4=(\mathcal{D}_4,\preceq)$ can be viewed by considering  the dotted edge (in red). }
     \label{Ex3}
 \end{figure}

Notice that the covering $\longrightarrow$ is a particular case of the Tamari covering
$\stackrel{\mathcal{T}}{\longrightarrow}$ whenever we take $Q=U^{k-1}D^{k-1}$ in the transformation $DUQD\stackrel{\mathcal{T}}{\longrightarrow} UQDD$ . Let $\mathcal{S}_n$ be the poset
$(\mathcal{D}_n,\leq)$.
The bottom part of Figure~\ref{fig1} shows an
example of such a covering, and Figure~\ref{Ex3} illustrates the Hasse
diagram of $\mathcal{S}_n$ for $n=4$ (without the red edges which
belongs to the Tamari lattice but not to this new poset). See also
Figure~\ref{fig6} for an illustration of the Hasse diagram of the Tamari
lattice and this poset for the case $n=6$.

In our study, the following facts will sometimes be used explicitly or implicitly.

\noindent{\bf Fact 1.} If $P\longrightarrow P'$,
then the path $P'$ is above the path $P$, that is, for any points
$(x,y)\in P'$ and $(x,z)\in P$ we have $y \geq z$.

\noindent{\bf Fact 2.}  If $P \longrightarrow P'$ then $\s(P)-1\leq\s(P') \leq \s(P)$ and $\st(P)\leq \st(P')+1$.

\section{Lattice structure of $\mathcal{S}_n=(\mathcal{D}_n,\leq)$.}
In this section,
we prove that the poset $\mathcal{S}_n=(\mathcal{D}_n,\leq)$ is a lattice and provide some related results.

\begin{lem} For $n\geq 2$, any Dyck path $P\in\mathcal{D}_n$, $P\neq U^nD^n$, contains at least one occurrence of $DU^kD^k$ for some $k\geq 1$.
  \label{le1}
\end{lem}
\noindent {\it Proof.} For $n\geq 2$, in any path from $\mathcal{D}_n$ not equal to $U^nD^n$, there exists an occurrence of $DU$, and the rightmost occurrence of $DU$ always starts an occurrence of $DU\, U^\ell D^{\ell} D$
  for some $\ell \geq 0$.
\hfill $\Box$

\begin{lem}
   For $n \geq 2$, any Dyck path $P\in\mathcal{D}_n$, $P\neq
  (UD)^n$, contains at least one occurrence of $U^kD^kD$ for some
   $k\geq 1$, and then $P$ contains at least one occurrence of $UDD$.
   \label{le2}
\end{lem}

\noindent {\it Proof.} By contradiction, let us assume that $P$
does not contain an occurrence $UDD$. This means that any peak $UD$ is
either at the end of $P$, or it precedes an up-step $U$, which implies
that a down-step cannot be contiguous to another down-step. Thus,
$P=(UD)^n$ which contradicts the hypothesis $P\neq (UD)^n$.  \hfill
$\Box$

\begin{prop}
  For any Dyck path $P\in\mathcal{D}_n$, we have $P\leq U^nD^n$ and
  $(UD)^n\leq P$, which means that the poset $(\mathcal{D}_n,\leq)$ admits a maximum element and a minimum element.
  \label{prop1}
\end{prop}

\noindent {\it Proof.} It suffices to apply Lemma~\ref{le1} and Lemma~\ref{le2}.
Indeed, let $P$ be a Dyck path in $\mathcal{D}_n$, $P\neq
U^nD^n$. Using Lemma~\ref{le1}, $P$ contains at least one occurrence of
$DU^kD^k$, $k\geq 1$. Let $P_1$ be the Dyck path obtained from $P$
after applying the covering $DU^kD^k\rightarrow U^kD^kD$ on this
occurrence. Due to Fact~1, any point $(x,y)$ in $P$ is below the point
$(x,z)$ in $P_1$ (i.e., $y\leq z$). Iterating the process from $P_1$,
we construct a sequence of coverings $P\rightarrow P_1\rightarrow
\ldots \rightarrow P_r$, $r\geq 1$, of Dyck paths that necessarily
converges towards a Dyck path without occurrence of $DU^kD^k$  for some
$k\geq 1$, i.e., towards $U^nD^n$, which implies $P\leq U^nD^n$.

By applying a similar argument, we can prove easily the second inequality.
\hfill $\Box$

\begin{prop}
  We consider $P,Q\in \mathcal{D}_n$ that satisfy $P\leq Q$, $P\neq Q$, and
  such that $P=RDS$ and $Q=RUS'$ ($R$ is the maximal common prefix of
  $P$ and $Q$). Let $W$ be the Dyck path obtained from $P$ by applying
  the covering $P\longrightarrow W$ on the leftmost occurrence of
  $DU^kD^k$, $k\geq 1$, in $DS$, then we necessarily have $W\leq
  Q$.
  \label{prop2}
\end{prop}

\noindent {\it Proof.} Let $P_0=P\rightarrow P_1 \rightarrow \ldots
\rightarrow P_k=Q$ be a sequence of coverings from $P$ to $Q$. 
Let us
suppose that the first covering is not applied on the leftmost
occurrence of $DU^kD^k$, $k\geq 1$, in $DS$. We call $D_0$ the first
step of this occurrence. Let $P_i\rightarrow P_{i+1}$ be the covering
involving $D_0$, and we write $P=R'D_0U^kD^kS'$ where $R$ is a prefix of $R'$ and $S'$ is a suffix of $S$. Necessarily, the down-step $D_0$ is followed by
$U^\ell D^\ell$ in $P_i$ where $\ell\geq k$. Then, all coverings
between $P$ and $P_i$ occur on the right of $D_0U^k$. We call them
$\alpha_1, \ldots ,\alpha_a$, $a\geq 1$. So, there is another sequence
of coverings from $P$ to $P_{i+1}$ by applying the following: we first
apply the covering $\beta$ involving $D_0U^kD^k$ in $P$, then we
obtain the path $R'U^kD^kDS'$; we apply the coverings $\alpha'_1,\ldots,
\alpha'_a$ (in the same order) where $\alpha'_i$ is the covering
involving the same occurrence of $DU^bD^b$ moved by $\alpha_i$, then
we obtain $R'U^kDU^{\ell-k}D^\ell S''$ for some $S''$; finally, we apply an additional
covering in order to obtain $P_{i+1}=R'U^\ell D^\ell DS''$ (see below
for an illustration of this process).

\begin{figure}[H]
\begin{center}
\begin{tikzpicture}[ultra thick]
 \draw (0,3) node {$P=R'D_0U^kD^k S'$};
 \draw [->] (2,3) -- (3,3);
 \draw [dashed] (3.5,3)--(4.5,3);
 \draw [->] (5,3) -- (6,3);
 \draw (8,3) node {$P_i=R'D_0U^\ell D^\ell S''$};
 \draw [->] (10,3) -- (11,3);
 \draw (13,3) node {$P_{i+1}=R'U^\ell D^\ell D S''$};
 \draw [->] (1,2.5) -- (1,1.3);
 \draw (0.5,0.5) node {$R'U^kD^kD S'$};
 \draw [->] (2,0.5) -- (3,0.5);
 \draw [dashed] (3.5,0.5)--(4.5,0.5);
 \draw [->] (5,0.5) -- (6,0.5);
 \draw (8,0.5) node {$R'U^kDU^{\ell-k}D^\ell  S''$};
 \draw [->] (9.7,0.8) -- (11.3,2.5);
 \draw (2.5,3.3) node {$\alpha_1$};
 \draw (5.5,3.3) node {$\alpha_a$};
 \draw (0.7,1.9) node {$\beta$};
 \draw (2.5,0.8) node {$\alpha'_1$};
 \draw (5.5,0.8) node {$\alpha'_a$};
 \end{tikzpicture}
 \end{center}
 \label{figcommute}
 \end{figure}

 Therefore, we can reach $Q$ from $P$ by first applying the covering involving $D_0U^kD^k$, which completes the proof.
\hfill $\Box$

As a byproduct of Proposition~\ref{prop2} and by a straightforward
induction, we have the following corollary.

\begin{cor}  The longest chain between $(UD)^n$ and $U^nD^n$ is of length $n(n-1)/2$.
\end{cor}
\noindent {\it Proof.} Thanks to Proposition~\ref{prop2}, the longest chain between $(UD)^n$ and $U^nD^n$ is unique and it can be constructed from $(UD)^n$ by applying at each step the leftmost covering. So, the length of this chain is given by $1+2+\ldots +(n-1)$, which gives the expected result.
\hfill $\Box$
\begin{thm} The poset $(\mathcal{D}_n,\leq)$ is a lattice.
  \label{tla}
\end{thm}
\noindent {\it Proof.}  For any $P,Q\in \mathcal{D}$, we need to prove that $P$ and $Q$ admit a join
and a meet element. Let us start to give the proof of the existence
of a join element. We proceed by induction on the semilength of the
Dyck paths. For $n\le3$,
  $\mathcal{S}_n=(\mathcal{D}_n,\leq)$ is isomorphic to the Tamari lattice.

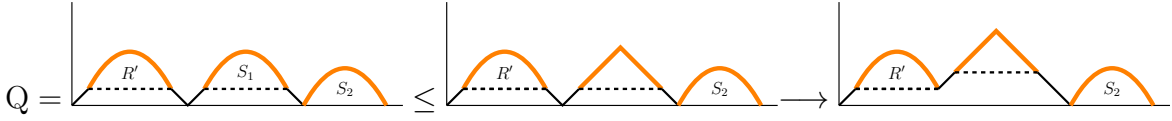
\begin{figure}[ht]
\begin{center}Q\;=\;\scalebox{0.55}{\begin{tikzpicture}[ultra thick]
 \draw[black, thick] (0,0)--(8,0);
 \draw[black, thick] (0,0)--(0,2.5);
 \draw[black, line width=1.5pt] (0,0)--(0.4,0.4); \draw[orange,very thick, line width=3pt] (0.4,0.4) parabola bend (1.4,1.3) (2.4,0.4);
  \draw[black, line width=1.5pt] (2.4,0.4)--(2.8,0)--(3.2,0.4);
  \draw[black, dashed] (0.4,0.4)--(2.4,0.4);
  \draw[orange,very thick, line width=3pt] (3.2,0.4) parabola bend (4.2,1.3) (5.2,0.4);
  \draw[black, dashed] (3.2,0.4)--(5.2,0.4);
  \draw[black, line width=1.5pt] (5.2,0.4)--(5.6,0);
  \draw[orange,very thick, line width=3pt] (5.6,0) parabola bend (6.6,0.9) (7.6,0.0);
   \draw  (1.4,0.8) node { $R'$};
   \draw  (4.2,0.8) node { $S_1$};
    \draw  (6.6,0.4) node { $S_2$};
 \end{tikzpicture}}\;${\leq}$\;\scalebox{0.55}{\begin{tikzpicture}[ultra thick]
 \draw[black, thick] (0,0)--(8,0);
 \draw[black, thick] (0,0)--(0,2.5);
 \draw[black, line width=1.5pt] (0,0)--(0.4,0.4); \draw[orange,very thick, line width=3pt] (0.4,0.4) parabola bend (1.4,1.3) (2.4,0.4);
  \draw[black, line width=1.5pt] (2.4,0.4)--(2.8,0)--(3.2,0.4);
  \draw[black, dashed] (0.4,0.4)--(2.4,0.4);
  %\draw[orange,very thick, line width=3pt] (3.2,0.4) parabola bend (4.2,1.3) (5.2,0.4);
  \draw[black, dashed] (0.4,0.4)--(2.4,0.4);
  \draw[orange, line width=3pt] (3.2,0.4)--(4.2,1.4)--(5.2,0.4);
  \draw[black, dashed] (3.2,0.4)--(5.2,0.4);
  \draw[black, line width=1.5pt] (5.2,0.4)--(5.6,0);
  \draw[orange,very thick, line width=3pt] (5.6,0) parabola bend (6.6,0.9) (7.6,0.0);
  \draw  (1.4,0.8) node { $R'$};
   %\draw  (4.2,0.8) node { $S_1$};
    \draw  (6.6,0.4) node { $S_2$};
 \end{tikzpicture}}$\longrightarrow$\;\scalebox{0.55}{\begin{tikzpicture}[ultra thick]
 \draw[black, thick] (0,0)--(8,0);
 \draw[black, thick] (0,0)--(0,2.5);
 \draw[black, line width=1.5pt] (0,0)--(0.4,0.4); \draw[orange,very thick, line width=3pt] (0.4,0.4) parabola bend (1.4,1.3) (2.4,0.4);
  \draw[black, line width=1.5pt] (2.4,0.4)--(2.8,0.8);
  \draw[black, dashed] (0.4,0.4)--(2.4,0.4);
  %\draw[orange,very thick, line width=3pt] (3.2,0.4) parabola bend (4.2,1.3) (5.2,0.4);
  \draw[black, dashed] (0.4,0.4)--(2.4,0.4);
  \draw[orange, line width=3pt] (2.8,0.8)--(3.8,1.8)--(4.8,0.8);
  \draw[black, dashed] (2.8,0.8)--(4.8,0.8);
  \draw[black, line width=1.5pt] (4.8,0.8)--(5.2,0.4)--(5.6,0);
  \draw[orange,very thick, line width=3pt] (5.6,0) parabola bend (6.6,0.9) (7.6,0.0);
  \draw  (1.4,0.8) node { $R'$};
  % \draw  (4.2,0.8) node { $S_1$};
    \draw  (6.6,0.4) node { $S_2$};
 \end{tikzpicture}}
\end{center}
\caption{An illustration of the construction $Q=UR'DUS_1DS_2\leq
  UR'DU^kD^kS_2\rightarrow UR'U^kD^kDS_2$ of the case ($ii$) in the
  proof of Theorem~\ref{tla}. }
\label{figth1_tmp}
\end{figure}

Now, let us assume that $\mathcal{S}_n=(\mathcal{D}_n,\leq)$ is a
lattice for $n\leq N$, and let us prove the result for $N+1$. Let $P$
and $Q$ be two paths in $\mathcal{D}_{N+1}$. We distinguish two cases
according to the form of the first return decompositions of $P$ and
$Q$.
\begin{enumerate}
    \item[(1)] If $P=URD S $ and $Q= UR'D S' $ where $R$ and $R'$ have
      the same length. Then we apply the induction hypothesis for $R$
      and $R'$ (resp. $S$ and $S'$), which means that $R\vee R'$
      (resp., $S\vee S'$) exists. Therefore, the path $ U (R\vee R') D
      (S\vee S') $ is necessarily the least upper bound of $P$ and
      $Q$, which proves that $P\vee Q$ exists.
    \item[(2)] Now, let us suppose that $P=URD S $ and $Q= UR'DS'$
      where the length $r'$ of $R'$ is strictly less than the length
      $r$ of $R$. Let $M$ be an upper bound of $P$ and $Q$ (there is
      at least one thanks to Proposition~\ref{prop1}). Since $r'<r$
      and due to Fact~1, $M$ has necessarily a decomposition
      $M=UM_1DM_2 $ where the length of $M_1$ is at least
      $r$.
      Therefore, in any sequence of coverings $Q\rightarrow
      \ldots \rightarrow M$ from $Q$ to $M$, there is necessarily a
      covering that involves and elevates the down-step just after
      $R'$. Due to the definition of the covering $\longrightarrow$,
      we can apply such a transformation only when this down-step is
      followed by an occurrence $U^kD^k$ for some $k\geq 1$.

    Assuming $S'=US_1DS_2$ and using Proposition~\ref{prop2}, we
    deduce that the inequality $Q=UR'DUS_1DS_2\leq M$ is equivalent to
    $UR'DU^kD^kS_2\leq M$ where $k\geq 1$ is the semilength of $US_1D$. Moreover, this condition is
    equivalent to $Q_1:=UR'U^kD^kDS_2\leq M$ (see
    Figure~\ref{figth1_tmp}). It is worth noting that $Q_1$ does not
    depend on the upper bound $M$. Iterating this process with $P$ and
    $Q_1$, we can construct two Dyck paths $P'$ and $Q'$ such that the
    condition $P\leq M$ and $Q\leq M$ is equivalent to $P'\leq M$ and
    $Q'\leq M$ where $P'$ and $Q'$ (that do not depend on $M$) are two
    Dyck paths lying in the first case of the proof.  Using the
    induction hypothesis, we conclude that $P'\vee Q'=P\vee Q$
    exists.
\end{enumerate}
 Considering the two cases, the induction is completed.

The existence of greatest lower bound then follows automatically since the poset is finite with a least element and a greatest element.
\hfill $\Box$

\section{Coverings, join and meet irreducible  elements}

For a given Dyck path $P$, the number  of incoming edges (in the Hasse
diagram) corresponds to the number $\st(P)$ of occurrences $U^aD^aD$, $a\geq
1$, in $P$, and the number of outgoing edges (coverings) corresponds
to the number $\s(P)$ of occurrences $DU^aD^a$, $a\geq 1$, in $P$.  Let
$A(x,y,z)$ be the trivariate generating function where the coefficient
of $x^ny^kz^\ell$ is the number of Dyck paths of semilength $n$ having
$k$ possible coverings (or equivalently $k$ outgoing edges), and $\ell$
incoming edges.

\begin{thm} We have
$$A(x,y,z)={\frac {R(x,y,z)-\sqrt {4\,x
 \left( xzy-xy-xz+1 \right)  \left( xy+xz-x-1 \right) +R(x,y,z) ^{2}}}{2x \left( xzy-xy
-xz+1 \right) }}
  ,$$
where $R(x,y,z)={x}^{2}zy-{x}^{2}y-{x}^{2}z+{x}^{2}-xy-xz+x+1$.
\label{th2}
\end{thm}
\noindent {\it Proof.} For short, we set $A:=A(x,y,z)$. We consider
the last return decomposition of a nonempty Dyck path $P$, that is
$P=R USD$ where $R$ and $S$ are two Dyck paths. We distinguish six
cases.
\begin{itemize}
    \item[(1)] If $R$ and $S$ are empty, then $P=UD$ and the g.f. for this path is $x$.
    \item[(2)] If $R$ is not empty and $S$ is empty, the g.f. for these paths is $(A-1)xy$.
    \item[(3)] If $R$ is empty and $S=U^aD^a$ with $a\geq 1$, then the g.f. for these paths is $\frac{x^2z}{1-xz}$.
    \item[(4)] If $R$ is not empty and $S=U^aD^a$ with $a\geq 1$, then the g.f. for these paths is $\frac{x^2z}{1-xz}(A-1)y$.
    \item[(5)] If $S=S'U^aD^a$ with $a\geq 1$ and $S'$ not empty ($R$ is possibly empty), then the g.f. for these paths is $\frac{x^2z}{1-xz}(A-1)yA$.
    \item[(6)] If $S$ does not end with $U^aD^a$, $a\geq 1$, then the g.f. for these paths is $AxB$ where $B$ is the g.f. for nonempty Dyck paths that do not end with a pyramid $U^aD^a$, $a\geq 1$. Using the complement, we have $B=A-1-x-\frac{x^2z}{1-xz}-x(A-1)y-\frac{x^2z}{1-xz}(A-1)y$ where $1$ is for the empty path; $x$ is for $UD$;
    $\frac{x^2z}{1-xz}$ is for the paths of the form $U^aD^a$, $a\geq 2$; $x(A-1)y$ is for the paths $RUD$ with $R$ not empty; and $\frac{x^2z}{1-xz}(A-1)y$ is for the paths $RU^aD^a$, $a\geq 2$.

    Summarizing the three cases, we obtain the functional equation
    \begin{align*}
        A=1+x+(A-1)xy+\frac{x^2z}{1-xz}+\frac{x^2z}{1-xz}(A-1)y+\frac{x^2z}{1-xz}(A-1)yA+\\
        +Ax\left(A-1-x-\frac{x^2z}{1-xz}-x(A-1)y-\frac{x^2z}{1-xz}(A-1)y\right),
    \end{align*}
    and a simple calculation provides the result.
\end{itemize}
\hfill $\Box$

The first terms of the series expansion are \begin{align*}1+x+ \left( z+y \right) {x}^{2}&+ \left( {y}^{2}+3\,yz+{z}^{2}
 \right) {x}^{3}+ \left( {y}^{3}+5\,{y}^{2}z+5\,y{z}^{2}+{z}^{3}+2\,yz
 \right) {x}^{4}+ \\
 &+\left( {y}^{4}+7\,{y}^{3}z+13\,{y}^{2}{z}^{2}+7\,y{z
}^{3}+{z}^{4}+5\,{y}^{2}z+5\,y{z}^{2}+3\,yz \right) {x}^{5}+O \left( {x}^{6} \right).
\end{align*}

\begin{rem}
  Observe that $R(x,y,z)=R(x,z,y)$ and $A(x,y,z)=A(x,z,y)$, which
  means that $(\st,\s)$ and $(\s,\st)$ are equidistributed on
  $\mathcal{D}_n$, $n\geq 1$.  In the following, we will show the existence of an
  involution on $\mathcal{D}_n$ that transports the bistatistic
  $(\st,\s)$ to $(\s,\st)$ (see Definition~\ref{def3}).
\end{rem}

\begin{cor} The generating function $E(x)$
where the coefficient of $x^n$ is the total number of possible
coverings over all Dyck paths of semilength $n$ (or equivalently the
number of edges in the Hasse diagram) is
  $$E(x)=\frac{-1+4 x +\left(1-2 x \right) \sqrt{1-4 x }}{2 \left(1-4 x \right) \left(1-x \right)}.$$
  The coefficient of $x^n$ is given by $$\sum_{k=0}^{n-2} \binom{2k+2}{k}.$$
  The ratio between the numbers of coverings in $\mathcal{T}_n$ and $\mathcal{S}_n$ tends towards $3/2$.
\end{cor}
\noindent {\it Proof.} We obtain $E(x)$ by calculating $\partial_y(A(x,y,1))\rvert_{y=1}$. Now we set $F(x)=\frac{E(x)(1-x)}{x}$.
Noticing that $F(x)=x\cdot \partial_x(C(x))$ where $C(x)$ is the generating function for the Catalan numbers satisfying $C(x)=1+xC(x)^2$. So, we deduce directly that
$$f_n:=[x^n]F(x)=\frac{n}{n+1}\binom{2n}{ n }.$$
Then, we have
$$[x^n]E(x)=\sum_{k=1}^{n-1}f_{k}=\sum_{k=1}^{n-1}\frac{k}{k+1}\binom{2k}{ k}=\sum_{k=0}^{n-2} \binom{2k+2}{k}.$$
Considering the asymptotics  of $[x^n]E(x)$ and  $(n-1) c_n/2$ (using classical methods, see~\cite{Fla} for instance), the limit of the ratio between the number of coverings in $\mathcal{T}_n$ and $\mathcal{S}_n$ is $3/2$.
\hfill $\Box$

 The first terms of the series expansion of $E(x)$ are $${x}^{2}+5\,{x}^{3}+20\,{x}^{4}+76\,{x}^{5}+286\,{x}^{6}+1078\,{x}^{7}
+4081\,{x}^{8}+15521\,{x}^{9}+O \left( {x}^{10} \right),$$ and the sequence of coefficients corresponds to \oeis{A057552} in~\cite{Sloa}.

Let $K(x)$ be the generating function where the coefficient of $x^n$ is the number of meet irreducible elements (Dyck paths with only one
 outgoing edge).  Using the symmetry $y\longleftrightarrow z$ in $A(x,y,z)$, $K(x)$ is also the generating function where the coefficient of $x^n$ is the number of join irreducible elements (Dyck paths with only one
 incoming edge).
\begin{cor} We have  $$K(x)=\frac{x^2}{(x-1)(2x-1)},$$ and the coefficient of $x^n$ is $2^{n-1}-1$ for $n\geq 1$.
\end{cor}
\noindent {\it Proof.} The generating function corresponds to the  coefficient of $y$ in the series expansion of $A(x,y,1)$, i.e., $[y]A(x,y,1)=\partial_y(A(x,y,1))\rvert_{y=0}$.
\hfill $\Box$

%Since the numbers of meet and join irreducible elements are the same, we deduce the following corollary.
%\begin{cor} The lattice $(\mathcal{D}_n,\leq)$ is extremal.
%    \end{cor}

Denote by $\mathcal{L}$ the set of paths with only one outgoing edge
and only one incoming edge.  Let $L(x)$ be the generating function
where the coefficient of $x^n$ is the number of such paths.
% being both meet and join irreducible elements

\begin{cor} We have  $$L(x)=3x^3+{\frac { \left( x+2 \right) {x}^{4}}{1-x-{x}^{2}}},$$
  and the coefficient of $x^n$ is $0$ whenever $n\leq 2$, is $3$
  whenever $n=3$, and is the Fibonacci number $F_{n-1}$ otherwise, where $F_n$
  is defined by $F_n=F_{n-1}+F_{n-2}$ with $F_1=F_2=1$.
\end{cor}
\noindent {\it Proof.} For $n\leq 2$, there is no such paths. For
$n=3$, there are three such paths, $UUDDUD$, $UDUUDD$ and
$UUDUDD$. For $n\geq 4$, the number of these paths is exactly the
coefficient of $y^2x^n$ in $A(x,y,y)$, {\it i.e.},
$[y^2x^n]A(x,y,y)={\frac { \left( x+2 \right) {x}^{4}}{1-x-{x}^{2}}}$.
\hfill $\Box$

\medskip

 For $n \ge 4$, we remark that any path from $\mathcal{L}$
  should contain exactly one occurrence of $DU^kD^k$, with $k\ge 1$, to
  guarantee the existence of a unique outgoing edge, and exactly one
  occurrence of $U^kD^kD$ for the incoming edge.  Thus, it must
  avoid overlapping $UDU$ anywhere except for the tail, which
  necessarily has a shape $U \, UD \, UD \, D^\ell$ for some $\ell > 0$.
  Figure~\ref{figfibo} presents a bijection between words of length
  $n\ge4$ from $\mathcal{L}$ and Knuth-Fibonacci words of length
  $n-3$, i.e. binary words avoiding consecutive 1s, discussed for
  example in Knuth's book~\cite[p. 286]{knuth3}. To construct the
  corresponding binary word, we read the Dyck path from left to right,
  write 1 for any up-step which starts a $UDU$ pattern, and 0
  otherwise.  The resulting word will always end with $011$, so we
  forget about these last 3 symbols.

      \begin{figure}[H]
        \begin{center}
          \scalebox{0.55}{
            \begin{tikzpicture}[ultra thick,
                l/.style={black, thick},
                b/.style={black, line width=1.5pt},
                d/.style={black, dashed, very thick}
              ]
              \draw[l] (0,0)--(13,0);

              \draw[b] (0,0)--(0.4,0.4);
              \draw[d] (0.4,0.4) -- (0.8,0.8);
              \draw[b] (0.8, 0.8)--(1.2,1.2)
              --(1.6, 0.8)--(2,1.2);
              \draw[d] (2,1.2) -- (2.8,2);
              \draw[b] (2.8, 2) -- (3.2, 2.4)
              -- (3.6, 2) -- (4, 2.4);
              \draw[d] (4,2.4) -- (4.8,3.2) --
              (5.2, 2.8) -- (5.6, 3.2);
              \draw[b] (5.6, 3.2) -- (6.4, 4)
              -- (6.8, 3.6) -- (7.2, 4) -- ($(7.2, 4) + 10*(+0.4, -0.4)$);

              \node at (0.1, 0.5) {\Large 0};
              \node at (0.4, 0.9) {\Large 0};
              \node at (0.7, 1.1) {\large $\cdot^{\cdot^\cdot}$};
              \node at (0.9, 1.3) {\Large 1};

              \node at (1.7, 1.3) {\Large 0};
              \node at (2.1, 1.7) {\large $\cdot^{\cdot^\cdot}$};
              \node at (2.5, 2.1) {\Large 0};
              \node at (2.9, 2.5) {\Large 1};

              \node at (3.7, 2.5) {\Large 0};
              \node at (4.1, 2.8) {\Large $\cdot^{\cdot^\cdot}$};

              \node at (5.7, 3.7) {\Large 0};
              \node at (6.1, 4.1) {\Large 1};
              \node at (6.9, 4.1) {\Large 1};

 \end{tikzpicture}}
        \end{center}
        \caption{Bijection between $\mathcal{L}$ and binary words
        avoiding consecutive 1s except the tail 011.}
        \label{figfibo}
 \end{figure}

\bigskip

In order to exhibit an involution on $\mathcal{D}_n$ that transports
the bistatistic $(\s,\st)$ to $(\st,\s)$, we need to define the
following subsets of $\mathcal{D}_n$. Let $\mathcal{D}_n^1$ (resp.,
$\mathcal{D}_n^2$), be the set of Dyck paths $P\in \mathcal{D}_n$ such
that $P$ ends with $DUD^k$, $k\geq 2$ (resp., ends with $DU^rD^k$,
$r\geq 2$, $k\geq 3$ and $r\neq k$). Less formally, $\mathcal{D}_n^1$ is the set of paths where  the
  last run of down-steps is of length at least two, and  the last run of up-steps is of length $1$; $\mathcal{D}_n^2$ is the set of paths where the
  last run of down-steps is of length at least two,  the last run of up-steps is of length at least $2$, and they do not end  with $U^\ell D^\ell$ for any $\ell \ge 1$. We set
$\mathcal{D}^i=\bigcup_{n\geq 0}\mathcal{D}^i_n$ for $i\in\{1,2\}$. For
instance, we have $UUDUDD\in \mathcal{D}_3^1$ and $UUUDDUUDDD\in
\mathcal{D}_5^2$.
%\marginpar{\color{red}\scriptsize Seems that M$\texttt{u}$\\
%okokok  is never used}
%{\color{red}\sout{For a Dyck path $P\in\mathcal{D}_n$, %we consider the
%statistic $\texttt{u}$ where $\texttt{u}(P)$ is the %length of its last
%run of up-steps. Obviously, we have $\texttt{u}(P)=1$ %whenever $P\in
%\mathcal{D}_n^1$, and $\texttt{u}(P)\geq 2$ whenever %$P\in
%\mathcal{D}_n^2$.}}

Below, we define recursively an involution $\phi$ from  $\mathcal{D}$ into itself as follows:
$$\left\{\begin{array}{lllr}
\phi(\epsilon)&= \epsilon,&&\hfill (i)\\
\phi(U^kD^k)&= (UD)^{k},& k\geq 1&\hfill (ii)\\
\phi(RU^kD^k)&= \phi(R)^\sharp (UD)^{k-1},&k\geq 2&\hfill (iii)\\
\phi(R_0UR_1\ldots UR_kUD^{k+1})&=\phi(R_0)U\phi(R_1)\ldots U\phi(R_k)UD^{k+1},& k\geq 1&\hfill  (iv)\\
\phi(R_0UR_1\ldots UR_kU^{r+1}D^{k+r+1})&=\phi(R_0)U\phi(R_1)\ldots U\phi(R_k)UD^{k+1}(UD)^r, &k,r\geq 1&\hfill  (v)\\
\end{array}\right.$$
where $R$ and $R_k$ are nonempty Dyck paths, and $R_i$, $0\leq i\leq k-1$, are Dyck paths (possibly empty), and $\phi(R)^\sharp$ is obtained from $\phi(R)$ by inserting an occurrence $UD$ on the last peak of $\phi(R)$, i.e., if $\phi(R)=R'UD^k$, then $\phi(R)^\sharp=R'UUDD^k$ (see below for the image by $\phi$ of paths of the form $RUD$).

    It is worth noticing that the definition of $\phi$ is given on a subset $\mathcal{P}_1$ of $\mathcal{D}$, consisting of Dyck paths of the forms $\epsilon$, $U^kD^k$ with $k\geq 1$, $RU^kD^k$ with $k\geq 2$ and $R$ nonempty, and Dyck paths in $\mathcal{D}^1\cup\mathcal{D}^2$. This means that the above definition of $\phi$ does not explicitly provide the image of paths of the form $RUD$, i.e. paths in $\mathcal{D}\backslash \mathcal{P}_1$.  However, this definition implicitly gives the image of paths in $\mathcal{D}\backslash \mathcal{P}_1$ by considering  the reverse of cases ($iii$) and ($v$) (this is possible because we have $\phi(\mathcal{D}^1)=\mathcal{D}^1$, and  $\phi(\mathcal{P}_1\backslash\mathcal{D}^1)=\mathcal{D}\backslash \mathcal{D}^1$). Then we obtain an involution $\phi$ on the entire set  $\mathcal{D}$. For instance, we have $\phi(UUUDDD)=UDUDUD$, $\phi(UUDDUD)=UDUUDD$, $\phi(UDUUDD)=UUDDUD$, $\phi(UUDUDD)=UUDUDD$ and $\phi(UDUDUD)=UUUDDD$.

\begin{figure}[h]
\begin{center}($ii$)\quad\scalebox{0.55}{\begin{tikzpicture}[ultra thick]
 \draw[black, thick] (0,0)--(11,0); \draw[black, thick] (0,0)--(0,2.5);

   \draw[orange, line width=3pt] (0,0)--(0.4,0.4)--(0.8,0.8)--(1.2,1.2)--(2.4,0);
   \draw  (2.8,0.8) node {$k\geq 1$};
 \end{tikzpicture}}\quad $\stackrel{\phi}{\longleftrightarrow}$\quad\scalebox{0.55}{\begin{tikzpicture}[ultra thick]
 \draw[black, thick] (0,0)--(11,0); \draw[black, thick] (0,0)--(0,2.5);
  %\draw[black, line width=3pt] (1.6,0)--(2,0.4);
   %\draw[orange, line width=3pt] (2,0.4)--(2.4,0.8)--(2.8,1.2)--(3.2,1.6)--(4.4,0.4);
  %\draw[black, line width=3pt] (4.4,0.4)--(4.8,0);
 %\draw[black, dashed, very thick] (2,0.4) -- (4.4,0.4);
\draw[orange, line width=3pt] (0,0)--(0.4,0.4)--(0.8,0)--(1.2,0.4)--(1.6,0)--(2,0.4)--(2.4,0);\draw  (1.5,0.8) node {$k\geq 1$};
 \end{tikzpicture}}
\end{center}
\begin{center}($iii$)\quad\scalebox{0.55}{\begin{tikzpicture}[ultra thick]
 \draw[black, thick] (0,0)--(11,0); \draw[black, thick] (0,0)--(0,2.5);
  \draw[orange,very thick, line width=3pt] (0,0) parabola bend (0.8,0.8) (1.6,0);
  %\draw[black, line width=1.5pt] (1.6,0)--(2,0.4);
   \draw[orange, line width=3pt] (1.6,0)--(2,0.4)--(2.4,0.8)--(2.8,1.2)--(3.2,1.6)--(4.4,0.4)--(4.8,0);
 % \draw[black, line width=1.5pt] (4.4,0.4)--(4.8,0);
 \draw[black, dashed, very thick] (2,0.4) -- (4.4,0.4);
 \draw  (0.8,1.2) node {$R$};\draw  (4.8,0.8) node {$k\geq 2$};
 \end{tikzpicture}}\quad $\stackrel{\phi}{\longleftrightarrow}$\quad\scalebox{0.55}{\begin{tikzpicture}[ultra thick]
 \draw[black, thick] (0,0)--(11,0); \draw[black, thick] (0,0)--(0,2.5);
  \draw[orange,very thick, line width=3pt] (0,0) parabola bend (0.8,0.8) (1.6,0);
  %\draw[black, line width=3pt] (1.6,0)--(2,0.4);
   %\draw[orange, line width=3pt] (2,0.4)--(2.4,0.8)--(2.8,1.2)--(3.2,1.6)--(4.4,0.4);
  %\draw[black, line width=3pt] (4.4,0.4)--(4.8,0);
 %\draw[black, dashed, very thick] (2,0.4) -- (4.4,0.4);
 \draw  (0.8,1.2) node {$\phi(R)^\sharp$};\draw[orange, line width=3pt] (1.6,0)--(2,0.4)--(2.4,0)--(2.8,0.4)--(3.2,0)--(3.6,0.4)--(4,0);\draw  (3,0.8) node {$k-1\geq 1$};
 \end{tikzpicture}}
\end{center}
\begin{center}($iv$)\quad\scalebox{0.55}{\begin{tikzpicture}[ultra thick]
 \draw[black, thick] (0,0)--(11,0); \draw[black, thick] (0,0)--(0,2.5);
  \draw[orange,very thick, line width=3pt] (0,0) parabola bend (0.8,0.8) (1.6,0);
  \draw[black, line width=1.5pt] (1.6,0)--(2,0.4);
  \draw[black, dashed, very thick] (2,0.4)--(3.3,0.4);
  \draw[orange,very thick, line width=3pt] (2,0.4) parabola bend (2.65,1) (3.3,0.4);
   \draw[black, line width=1.5pt] (3.3,0.4)--(3.7,0.8);
   \draw[black, dashed, very thick] (3.7,0.8) -- (4.9,0.8);
   \draw[black, line width=1.5pt] (4.9,0.8)--(5.3,1.2);
   \draw[orange,very thick, line width=3pt] (5.3,1.2) parabola bend (5.95,1.8) (6.6,1.2);
   \draw[black, dashed, very thick] (5.3,1.2)--(6.6,1.2);

   \draw[black, line width=1.5pt] (6.6,1.2)--(7,1.6)--(8.7,0);
 \draw  (0.9,1.5) node {$R_0$};\draw  (2.6,1.5) node { $R_1$};
 \draw  (5.85,2.4) node { $R_k$};
 \end{tikzpicture}}\quad $\stackrel{\phi}{\longleftrightarrow}$\quad\scalebox{0.55}{\begin{tikzpicture}[ultra thick]
 \draw[black, thick] (0,0)--(11,0); \draw[black, thick] (0,0)--(0,2.5);
  \draw[orange,very thick, line width=3pt] (0,0) parabola bend (0.8,0.8) (1.6,0);
  \draw[black, line width=1.5pt] (1.6,0)--(2,0.4);
  \draw[black, dashed, very thick] (2,0.4)--(3.3,0.4);
  \draw[orange,very thick, line width=3pt] (2,0.4) parabola bend (2.65,1) (3.3,0.4);
   \draw[black, line width=1.5pt] (3.3,0.4)--(3.7,0.8);
   \draw[black, dashed, very thick] (3.7,0.8) -- (4.9,0.8);
   \draw[black, line width=1.5pt] (4.9,0.8)--(5.3,1.2);
   \draw[orange,very thick, line width=3pt] (5.3,1.2) parabola bend (5.95,1.8) (6.6,1.2);
   \draw[black, dashed, very thick] (5.3,1.2)--(6.6,1.2);
   \draw[black, line width=1.5pt] (6.6,1.2)--(7,1.6)--(8.7,0);
 \draw  (0.9,1.5) node {$\phi(R_0)$};\draw  (2.6,1.5) node { $\phi(R_1)$};
 \draw  (5.85,2.4) node { $\phi(R_k)$};
 \end{tikzpicture}}
\end{center}
\begin{center}($v$)\quad\scalebox{0.55}{\begin{tikzpicture}[ultra thick]
 \draw[black, thick] (0,0)--(11,0); \draw[black, thick] (0,0)--(0,2.5);
  \draw[orange,very thick, line width=3pt] (0,0) parabola bend (0.8,0.8) (1.6,0);
  \draw[black, line width=1.5pt] (1.6,0)--(2,0.4);
  \draw[black, dashed, very thick] (2,0.4)--(3.3,0.4);
  \draw[orange,very thick, line width=3pt] (2,0.4) parabola bend (2.65,1) (3.3,0.4);
   \draw[black, line width=1.5pt] (3.3,0.4)--(3.7,0.8);
   \draw[black, dashed, very thick] (3.7,0.8) -- (4.9,0.8);
   \draw[black, line width=1.5pt] (4.9,0.8)--(5.3,1.2);
   \draw[orange,very thick, line width=3pt] (5.3,1.2) parabola bend (5.95,1.8) (6.6,1.2);
   \draw[black, dashed, very thick] (5.3,1.2)--(6.6,1.2);
   \draw[orange, line width=3pt] (7,1.6)--(7.8,2.4)--(8.6,1.6);
   \draw[black, dashed, very thick] (7,1.6)--(8.6,1.6);
   \draw[black, line width=1.5pt] (6.6,1.2)--(7,1.6);
   \draw[black, line width=1.5pt] (8.6,1.6)--(10.2,0);
    \draw  (0.9,1.5) node {$R_0$};\draw  (2.6,1.5) node { $R_1$};
   \draw  (5.85,2.4) node { $R_k$};
    \draw  (9,2) node { $r\geq 1$};
 \end{tikzpicture}}\quad $\stackrel{\phi}{\longleftrightarrow}$\quad\scalebox{0.55}{\begin{tikzpicture}[ultra thick]
 \draw[black, thick] (0,0)--(11,0); \draw[black, thick] (0,0)--(0,2.5);
  \draw[orange,very thick, line width=3pt] (0,0) parabola bend (0.8,0.8) (1.6,0);
  \draw[black, line width=1.5pt] (1.6,0)--(2,0.4);
  \draw[black, dashed, very thick] (2,0.4)--(3.3,0.4);
  \draw[orange,very thick, line width=3pt] (2,0.4) parabola bend (2.65,1) (3.3,0.4);
   \draw[black, line width=1.5pt] (3.3,0.4)--(3.7,0.8);
   \draw[black, dashed, very thick] (3.7,0.8) -- (4.9,0.8);
   \draw[black, line width=1.5pt] (4.9,0.8)--(5.3,1.2);
   \draw[orange,very thick, line width=3pt] (5.3,1.2) parabola bend (5.95,1.8) (6.6,1.2);
   \draw[black, dashed, very thick] (5.3,1.2)--(6.6,1.2);
   \draw[black, line width=1.5pt] (6.6,1.2)--(7,1.6)--(8.7,0);
    \draw[orange, line width=3pt] (8.7,0)--(9.1,0.4)--(9.5,0)--(9.9,0.4)--(10.3,0);
 \draw  (0.9,1.5) node {$\phi(R_0)$};\draw  (2.6,1.5) node { $\phi(R_1)$};\draw  (9.5,0.8) node { $r\geq 1$};
 \draw  (5.85,2.4) node { $\phi(R_k)$};
 \end{tikzpicture}}
\end{center}
\caption{An illustration of the involution $\phi$ for the cases ($ii$)-($v$).
}
\label{fig2}\end{figure}
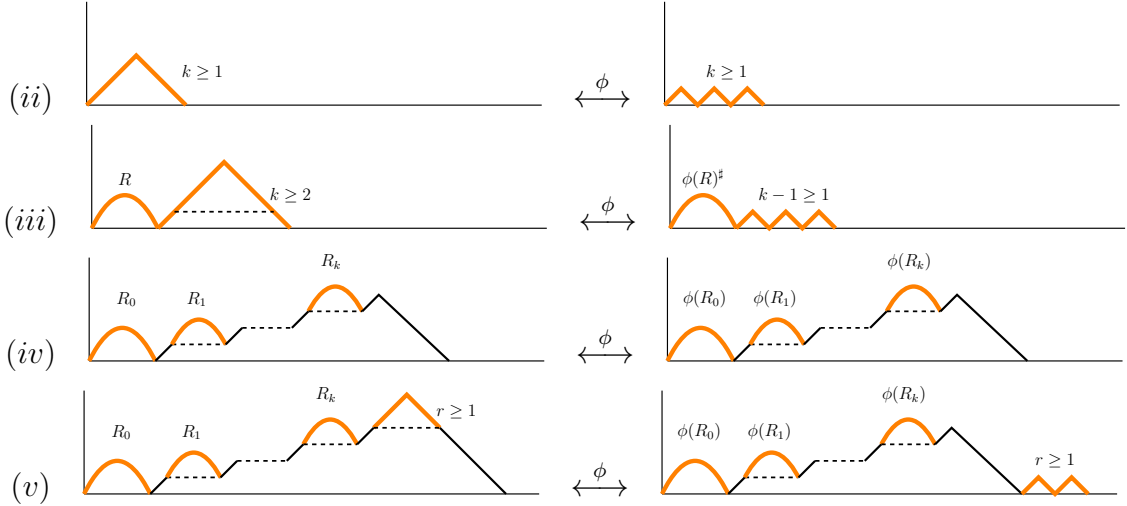

\bigskip

 Let $P$ be a Dyck path, we recall (see Definition~\ref{def3}) that $\texttt{s}(P)$ (resp. $\texttt{t}(p)$) gives the number of occurrences $DU^kD^k$ (resp. $U^kD^kD$), $k\geq 1$, in $P$.

\begin{thm} The map $\phi$ is an involution on $\mathcal{D}$ that preserves the semilength, and that transports the bistatistic $(\s,\st)$ to $(\st,\s)$. Moreover, $\phi$ preserves the number of asymmetric pyramids, and preserves the sum of the weights of symmetric pyramids.
\label{th3}
\end{thm}

\noindent {\it Proof.} We proceed by induction on the semilength $n$
of the paths. We assume the result for paths of semilength at most
$n$, and we prove the result for the semilength $n+1$.  We distinguish
four cases according to the definition of $\phi$ (we omit the case ($i$)
since $n+1\geq 1$).

\begin{enumerate}
    \item[($ii$)] If $P=U^{n+1}D^{n+1}$ then we have
      $\phi(P)=(UD)^{n+1}$, $\s(P)=0=\st(\phi(P))$, and
      $\st(P)=n=\s(\phi(P))$.
    \item[($iii$)] If $P=RU^{k}D^{k}$, with $R$ nonempty and $k\geq
      2$, then we have $\phi(P)=\phi(R)^\sharp(UD)^{k-1}$,
      $\s(P)=\st(R)+1$, and using the induction hypothesis, this is
      equal to $\st(\phi(R))+1=\st(\phi(P)^\sharp)$. Moreover, we
      have $\st(P)=\st(R)+k-1$ and using the induction hypothesis,
      this is equal to
      $\s(\phi(R))+k-1=\s(\phi(P)^\sharp)+k-1=\s(\phi(P)^\sharp(UD)^{k-1})=\s(\phi(P))$.
     \item[($iv$)] If $P=R_0UR_1\ldots UR_kUD^{k+1}$, then we have
       $\phi(P)=\phi(R_0)U\phi(R_1)\ldots U\phi(R_k)UD^{k+1}$, and
       $\s(P)=\sum_{i=0}^{k}\s(R_i)+1$, and using the induction 
       hypothesis, this is equal to
       $\sum_{i=0}^{k}\st(\phi(R_i))+1=\st(\phi(R_0)U\phi(R_1)\ldots
       U\phi(R_k)UD^{k+1})=\st(\phi(P))$. A similar argument allows us to
       prove $\st(P)=\s(\phi(P))$.
    \item[($v$)] If $P=R_0UR_1\ldots UR_kUD^{k+1}$, then we have
      $\phi(P)=\phi(R_0)U\phi(R_1)\ldots U\phi(R_k)UD^{k+1}(UD)^r$,
      and $\s(P)=\sum_{i=0}^{k}\s(R_i)+1$, and using the induction
      hypothesis, this is equal to
      $\sum_{i=0}^{k}\st(\phi(R_i))+1=\st(\phi(R_0)U\phi(R_1)\ldots
      U\phi(R_k)UD^{k+1}(UD)^r)=\st(\phi(P))$. The equality
      $\st(P)=\s(\phi(P))$ is obtained in the same way.
\end{enumerate} Considering all these cases, we obtain the result by induction.

\hfill $\Box$

Notice that Theorem~\ref{th3} allows us to retrieve the symmetry
obtained by Theorem~\ref{th2}, but unfortunately, the involution
$\phi$ does not induce a symmetry on the lattice
$(\mathcal{D}_n,\leq)$.

\section{Enumeration of intervals}
In this section, we provide the generating function and a closed form
for the number of intervals in the lattice $\mathcal{S}_n$. The method
is inspired by the work of Bousquet-M{\'e}lou and
Chapoton~\cite{Boucha}. We will use a catalytic variable that
considers the size of last run of down-steps.  We introduce the
bivariate generating function $$I(x,y)=\sum\limits_{n,k\geq 1}
a_{n,k}x^ny^k,$$ where $a_{n,k}$ is the number of
intervals in $\mathcal{S}_n$ such that the upper path ends with $k$
down-steps exactly. We also define $$J(x,y)=\sum\limits_{n,k\geq 1}
b_{n,k}x^ny^k,$$ where $b_{n,k}$ is the number of
intervals in $\mathcal{S}_n$ such that the upper path is prime and
ends with $k$ down-steps exactly (recall that a Dyck path is prime
whenever it only touches the $x$-axis at its beginning and its end).

\begin{lem} The following functional equation holds
    $$I(x,y)=J(x,y)+I(x,1)\cdot J(x,y).$$
    \label{le3}
\end{lem}
\noindent {\it Proof.} Let $(P,Q)$ be an interval in $\mathcal{S}_n$ where $P$ is the lower bound and $Q$ the upper bound.
We distinguish two cases.
\begin{enumerate}
    \item[(1)] If $Q$ is prime, then the contribution for these intervals is simply $J(x,y)$.
    \item[(2)] Otherwise, $Q$ is not prime and it has a last return decomposition $Q=R U S D$ with $R,S\in\mathcal{D}$ and $R$ not empty. This implies that $P$ is of the form $P=P_1P_2$ where $P_1,P_2\in\mathcal{D}$, and $P_2$ and $USD$ have the same length. This induces a bijection from intervals in this case with pairs of intervals $I_1:=(P_1,R)$ and $I_2:=(P_2,USD)$, where $I_1$ is any interval of of smaller length, and $I_2$ is any interval (of smaller length) lying in the case ($1$). Therefore, the contribution for the intervals in this case is $I(x,1)\cdot J(x,y)$.
\end{enumerate}
Considering the two cases, we obtain the expected result.
\hfill $\Box$

\begin{lem} The following functional equation holds
    $$J(x,y)=xy+xyI(x,y)+\frac{J(x,y)-J(x,1)}{y-1}\cdot C(xy)xy^2,$$
    where $C(x)$ is the g.f. for Catalan numbers, i.e., $C(x)=1+xC(x)^2$.
    \label{le4}
\end{lem}
\noindent {\it Proof.} Let $(P,Q)$ be an interval in $\mathcal{S}_n$ where $P$ is the lower bound and $Q$ the upper bound.
We distinguish three cases.
\begin{enumerate}
    \item[(1)] If $P=UD$ and $Q=UD$, then the contribution is $xy$.
    \item[(2)] If $P$ is prime, then we have $P=UP'D$ and we necessarily have $Q=UQ'D$ (i.e., $Q$ is prime) where $P'$ and $Q'$ are nonempty Dyck paths. Thus, the contribution for these intervals is $xyI(x,y)$.
    \item[(3)] Otherwise, $P$ is not prime which means that it has a last return decomposition $P=R U S D$ with $R,S\in\mathcal{D}$ and $R$ not empty (see the bottom of Figure~\ref{figformQ}). As $Q$ is prime, in any path of coverings $P\rightarrow P_1\rightarrow \cdots\rightarrow Q$ from $P$ to $Q$, there is necessarily a covering that involves and elevates the up-step just after $R$. We can apply such a transformation only when the suffix $P$ delimited by this up-step is of the form $U^kD^k$ for some $k\geq 1$.  This condition implies that $Q$ is necessarily of the following form
    $Q=Q'U^kD^{k+\ell}$ where $Q'$ is a prefix of Dyck path and $\ell\geq 1$ (see the top of Figure~\ref{figformQ} for an illustration of the form of $Q=Q'U^kD^{k+\ell}$). Let $h\geq 1$ be the height of the right point of the  last up-step of $Q'$. Then, we necessarily have $h\geq \ell\geq 1$.

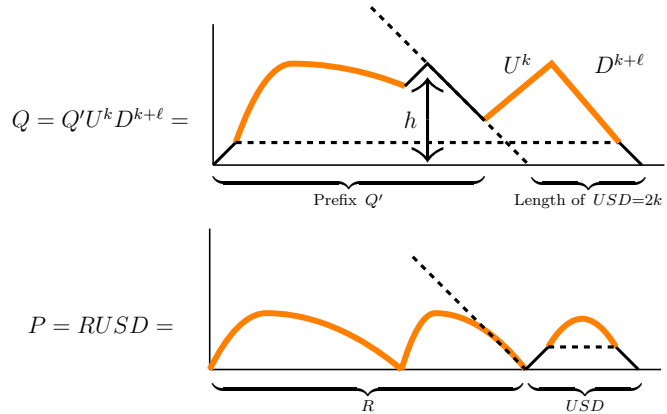
\begin{figure}[H]
\begin{center}
\scalebox{0.75}{\begin{tikzpicture}[ultra thick]
 \draw[black, thick] (0,0)--(8,0);
 \draw[black, thick] (0,0)--(0,2.5);
 \draw[black, line width=1.5pt] (0,0)--(0.4,0.4); \draw[orange,very thick, line width=3pt] (0.4,0.4) parabola bend (1.4,1.8) (3.4,1.4);
 % \draw[black, line width=3pt] (2.4,0.4)--(2.8,0.8);
  %\draw[black, dashed] (0.4,0.4)--(2.4,0.4);
  %\draw[orange,very thick, line width=3pt] (3.2,0.4) parabola bend (4.2,1.3) (5.2,0.4);
  \draw[black, dashed] (0.4,0.4)--(7.2,0.4);
  \draw[black, line width=1.5pt] (3.4,1.4)--(3.8,1.8)--(4.8,0.8);
  \draw[orange, line width=3pt] (4.8,0.8)--(6,1.8)--(7.2,0.4);
 % \draw[black, dashed] (2.8,0.8)--(4.8,0.8);
  \draw[black, dashed] (2.8,2.8)--(5.2,0.4)--(5.6,0);
  \draw[black, line width=1.5pt] (7.2,0.4)--(7.6,0);
  %\draw[orange,very thick, line width=3pt] (5.6,0) parabola bend (6.6,0.9) (7.6,0.0);
  \draw  (5.4,1.8) node {$U^k$};\draw  (7.2,1.8) node {$D^{k+\ell}$};
  \draw  (3.5,0.8) node {$h$};
   \draw  (6.65,-0.5) node {$\underbrace{\hskip2cm}_{\text{Length of}~USD=2k} $};
   \draw  (2.4,-0.5) node {$\underbrace{\hskip4.8cm}_{\text{Prefix}~Q'} $};
   \draw  (3.8,0.8) node {\scalebox{3}{$\big\updownarrow$}};
    \draw  (-2,0.8) node {$Q=Q'U^kD^{k+\ell}=$};
 \end{tikzpicture}}
\end{center}
\begin{center}
\scalebox{0.75}{\begin{tikzpicture}[ultra thick]
 \draw[black, thick] (0,0)--(8,0);
 \draw[black, thick] (0,0)--(0,2.5);
 %\draw[black, line width=3pt] (0,0)--(0.4,0.4);
 \draw[orange,very thick, line width=3pt] (0,0) parabola bend (1,1)  (3.4,0);
 \draw[orange,very thick, line width=3pt] (3.4,0) parabola bend (4,1)  (5.6,0);
  \draw[black, line width=1.5pt] (5.6,0)--(6,0.4);
  \draw[black, line width=1.5pt] (7.2,0.4)--(7.6,0);
   \draw[black, dashed] (3.6,2)--(5.2,0.4)--(5.6,0);
  \draw[orange,very thick, line width=3pt] (6,0.4) parabola bend (6.6,0.9) (7.2,0.4);
  \draw[black, dashed] (6,0.4)--(7.2,0.4);
   \draw  (6.65,-0.5) node {$\underbrace{\hskip2cm}_{USD} $};
   \draw  (2.8,-0.5) node {$\underbrace{\hskip5.5cm}_{R} $};
    \draw  (-2,0.8) node {$~~~P=RUSD=$};
 \end{tikzpicture}}
\end{center}
\caption{The form of the upper bound $Q=Q'U^kD^{k+\ell}$, $\ell\geq 1$, and the form of the lower bound $P=RUSD$. }
\label{figformQ}
\end{figure}

  On the other hand, any $Q$ of this form (i.e., $Q=Q'U^kD^{k+\ell}$ with $h\geq \ell\geq 1$) is candidate for an upper bound of an interval $(P,Q)=(RUSD,Q'U^kD^{k+\ell})$ if and only if $(R,Q'D^\ell)$ and $(USD, U^kD^k)$ are two intervals, which is equivalent to the condition $(R,Q'D^\ell)$ is an interval (indeed, $(USD, U^kD^k)$  is always an interval).

    So, for each interval of the form $(R,Q'D^\ell)$, $h\geq \ell \geq 1$, we can construct $h$ intervals of the form $(RUSD,Q''U^kD^{k+\ell}$, $1\leq \ell\leq h$, where $Q''$ is the greatest prefix of $Q'$ ending by an up-step, i.e., the intervals  $$(RUSD,Q''D^{h-1}U^kD^{k+1}),$$ $$(RUSD,Q''D^{h-2}U^kD^{k+2}),$$ $$\ldots \ldots$$ $$(RUSD,Q''DU^kD^{k+h-1}),$$ $$(RUSD,Q''U^kD^{k+h}).$$

    Considering this study, the contribution for the intervals in this case is given by

    $$xyC(xy)\cdot\sum_{h\geq 1}\sum_{n\geq 1}b_{n,h}x^n(y+y^2+\ldots +y^h)=
xy^2C(xy)\cdot\sum_{h\geq 1}\sum_{n\geq 1}b_{n,h}x^n\frac{y^h-1}{y-1},
$$
that can be expressed as
$$xy^2C(xy)\cdot\left(\sum_{h\geq 1}\sum_{n\geq 1}b_{n,h}x^n\frac{y^h}{y-1}-\sum_{h\geq 1}\sum_{n\geq 1}b_{n,h}x^n\frac{1}{y-1}\right)=
C(xy)xy^2\cdot\frac{J(x,y)-J(x,1)}{y-1}.
$$

\end{enumerate}
Considering the three cases, we obtain the expected result.
\hfill $\Box$

\bigskip

From Lemma~\ref{le3} and Lemma~\ref{le4} and after a straightforward substitution,  we obtain a system of equations.

\begin{thm} The following system of functional equations holds:
    $$\left\{\begin{array}{ll}
I(x,y)&=\frac{J(x,y)}{1-J(x,1)},\\
J(x,y)&=xy+xy\frac{J(x,y)}{1-J(x,1)}+\frac{J(x,y)-J(x,1)}{y-1}\cdot C(xy)xy^2.\\
\end{array}\right.$$
\end{thm}

In order to compute $J(x,1)$, we use the kernel method~\cite{ban, pro}
on $$J(x,y)\cdot
\left(1-\frac{xy}{1-J(x,1)}-\frac{C(xy)xy^2}{y-1}\right)=xy-\frac{J(x,1)}{y-1}\cdot
C(xy)xy^2.$$ This method consists in cancelling the factor of $J(x,y)$
by finding $y$ as an algebraic function $y_0$ of $J(x,1)$ and $x$. So,
if we substitute $y_0$ for $y$ in the right hand side of the
equation, then it necessarily equals zero (in order to counterbalance the
cancellation on the left-hand side). So, we deduce the following system
of equations that allows us to determine $J(x,1)$, and then $J(x,y)$:

$$\left\{\begin{array}{ll}
1-\frac{xy_0}{1-J(x,1)}-\frac{C(xy_0)xy_0^2}{y_0-1}&=0,\\
xy_0-\frac{J(x,1)}{y_0-1}\cdot C(xy_0)xy_0^2&=0.\\
\end{array}\right.$$

So, we deduce $y_0=\frac{1+4x-\sqrt{1-8x}}{8x}$, and the following theorem.

\begin{thm}The generating function $J(x,y)$ for the number of intervals $(P,Q)$ where $Q$ is prime,  with respect to the semilength and the size of the last run of down-steps is
\begin{align*}
  J(x,y)&={\frac {xy \left( -1+{\it J(x,1)} \right)  \left( {\it J(x,1)}\,{\it C}
 \left( xy \right) y-y+1 \right) }{{\it J(x,1)}\,{\it C} \left( xy
 \right) x{y}^{2}-{\it C} \left( xy \right) x{y}^{2}-x{y}^{2}-{\it J(x,1)
}\,y+xy+{\it J(x,1)}+y-1}}
%&=\frac{\left(3+\sqrt{1-8 x}\right) y \left(\left(\sqrt{1-8 x}-1\right) \left(-1+\sqrt{-4 x y +1}\right)-8 x y +8 x\right)}{-4 \left(3+\sqrt{1-8 x}\right) y \sqrt{-4 x y +1}+\left(-4 y +8\right) \sqrt{1-8 x}+24+32 y^{2} x +\left(-32 x -12\right) y}.\\
 \end{align*}
with $$J(x,1)=\frac{1-\sqrt{1-8x}}{4},$$
and $C(x)$ is the g.f. for Catalan numbers, i.e. $C(x)=1+xC(x)^2$.
\end{thm}

The series expansion of $J(x,1)$
is $$x+2\,{x}^{2}+8\,{x}^{3}+40\,{x}^{4}+224\,{x}^{5}+1344\,{x}^{6}+8448\,
{x}^{7}+54912\,{x}^{8}+366080\,{x}^{9}+O \left( {x}^ {10} \right),$$
where the sequence of coefficients corresponds to the sequence
\oeis{A052701} in~\cite{Sloa} that counts outerplanar maps with a
given number of edges~\cite{Geff}. The $n$-th coefficient is given by the closed form
$$2^{n-1}c_{n-1},$$ where $c_n=(2n)!/(n!(n+1)!)$ is the $n$-th Catalan number \oeis{A000108} in~\cite{Sloa}.

The series expansion of $J(x,y)$ is
\begin{align*} yx+2\,{y}^{2}{x}^{2}+ \left( 5\,y+3 \right) {y}^{2}{x}^{3}+ \left( 14
\,{y}^{2}+15\,y+11 \right) {y}^{2}{x}^{4}+ \left( 42\,{y}^{3}+61\,{y}^
{2}+68\,y+53 \right) {y}^{2}{x}^{5}+\\
\left( 132\,{y}^{4}+233\,{y}^{3}+
325\,{y}^{2}+363\,y+291 \right) {y}^{2}{x}^{6}+O \left( {x}^{7} \right),
\end{align*}

\begin{thm} The generating function $I(x,y)$ for the number of intervals $(P,Q)$  with respect to the semilength and the size of the last run of down-steps is
$$I(x,y)=J(x,y)\cdot{\frac { 3- \sqrt {1-8x}  }{2(x+1)}},$$
and the generating function $I(x,1)$ for the number of intervals $(P,Q)$  with respect to the semilength is
$$I(x,1)={\frac {1-2\,x-\sqrt {1-8x}}{2(x+1)}}.$$
\end{thm}

The series expansion of $I(x,1)$ is
$$x+3\,x^2+13\,x^3+67\,x^4+381\,x^5+2307\,x^6+14589
\,x^7+95235\,x^8+636925\,x^9+O(x^{10}),$$ where the sequence of
coefficients corresponds to the sequence \oeis{A064062} in~\cite{Sloa} that counts
simple outerplanar maps with a given number of vertices~\cite{Geff}. The $n$-th coefficient is given by the closed form $$\frac{1}{n}\sum\limits_{m=0}^{n-1}(n-m)\binom{n+m-1}{m}2^m.$$ An asymptotic approximation for the ratio of the numbers of intervals in $\mathcal{T}_n$ and $\mathcal{S}_n$ is
$$\frac{2^{5n+\frac{5}{2}}}{n\cdot 3^{3n+\frac{1}{2}}}.$$

The series expansion of $I(x,y)$ is
\begin{align*}
  yx+ \left( 2\,y+1 \right) yx^2+ \left( 5\,{y}^2+5\,y+3 \right)
yx^3+ \left( 14\,{y}^3+20\,{y}^2+20\,y+13 \right) yx^4+\\
 \left( 42\,{y}^4+75\,{y}^3+98\,{y}^2+99\,y+67 \right) yx^5+O \left( x^6
 \right).
\end{align*}
\section{Going further}

We discussed here several open questions related to the lattice $\mathcal{S}_n=(\mathcal{D}_n, \le)$.

\bigskip

\begin{conj} Finding a combinatorial interpretation of the equality $J(x,1)=x+2J(x,1)^2$.
\end{conj}

An {\it outerplanar map}~\cite{Geff} is a connected planar multigraph
with a specific embedding in the $2$-sphere, up to oriented
homeomorphisms, where a root edge is selected and oriented, and such
that all its vertices are in the outer face.

\begin{conj} Find a nice bijection between
  intervals in $\mathcal{S}_n$ and the simple outerplanar maps
  with a given number of vertices.
\end{conj}

The {\it distance} between two Dyck paths $P$ and $Q$ is the length of
a shortest path between $P$ in $Q$ in the underlying undirected graph
of the poset.

\begin{conj}
  Is there a polynomial time algorithm to compute the distance between
  two Dyck paths in $\mathcal{S}_n$?
\end{conj}

The {\it diameter} is the maximum distance between any two vertices.

\begin{conj} For $n\geq 3$,  we conjecture that the diameter of $\mathcal{S}_n$  is $2n-4$, and that this value corresponds to  the distance between $(UD)^n$ and    $UU(UD)^{n-2}DD$.
\end{conj}

The {\it M\"obius function}~\cite{Stamob} of $\mathcal{S}_n$, $\mu:\mathcal{S}_n\rightarrow \mathbb{Z}$ is defined recursively by
$$\mu(P)=-\sum\limits_{Q<P}\mu(Q)~\mbox{ if } ~P\neq (UD)^n,~\mbox{  with initial condition } \mu((UD)^n)=1.$$

\begin{conj} For $n\geq 2$, is there an efficient and non-recursive algorithm to compute the M\"obius function of $\mathcal{S}_n$ as in~\cite{palmob}?
\end{conj}

An $m$-Dyck path of size $n$ is a path on $\mathbb{N}^2$, starting at the origin and ending at $(2nm,0)$, consisting of up-steps
$(m,m)$ and down-steps $(1,-1)$. In the literature, these paths have been endowed with a lattice structure that generalizes the Tamari lattice~\cite{Berg,Bous}.

\begin{conj} Is it possible to generalize this study for $m$-Dyck paths?
\end{conj}

The Tamari covering transforms an occurrence of $DUPD$ into $UPDD$
without any constraints on the structure of a Dyck path $P$.  The
covering relation discussed in this paper can be viewed as the
following restriction of the Tamari covering: we only allow $DUPD
\longrightarrow UPDD$ whenever the path $P$ avoids consecutive pattern
$DU$. This opens the way to the further generalization:
{\em pattern-avoiding Tamari poset} is given by the transitive closure
of the covering $DUPD \longrightarrow UPDD$ where $P$ avoids a given
consecutive pattern $\nu$. It is easy to prove that any Dyck path of
semilength $n$ can be obtained from $(UD)^n$ by a sequence of such
pattern-aware transformations for any $\nu$.  For $k \ge 1, n \ge
k+2$, if $\nu = U^k$, then the poset will have at least two non
comparable maximal elements $U^nD^n$ and $UD U^{k+1}D^{k+1}$.  In
general the situation is not that simple, for example, if $\nu = UDU$
the poset has a minimum and maximum element, but some pairs of paths
do not have a greatest lower bound. However, if $\nu=UUDU$ then it seems that we obtain a
lattice structure. So, it becomes natural to ask the following.

\begin{conj} Could we characterize patterns $\nu$  inducing a lattice structure in the pattern-avoiding  version of the Tamari poset?
\end{conj}

\begin{figure}[ht]
     \centering \includegraphics[scale = 0.2]{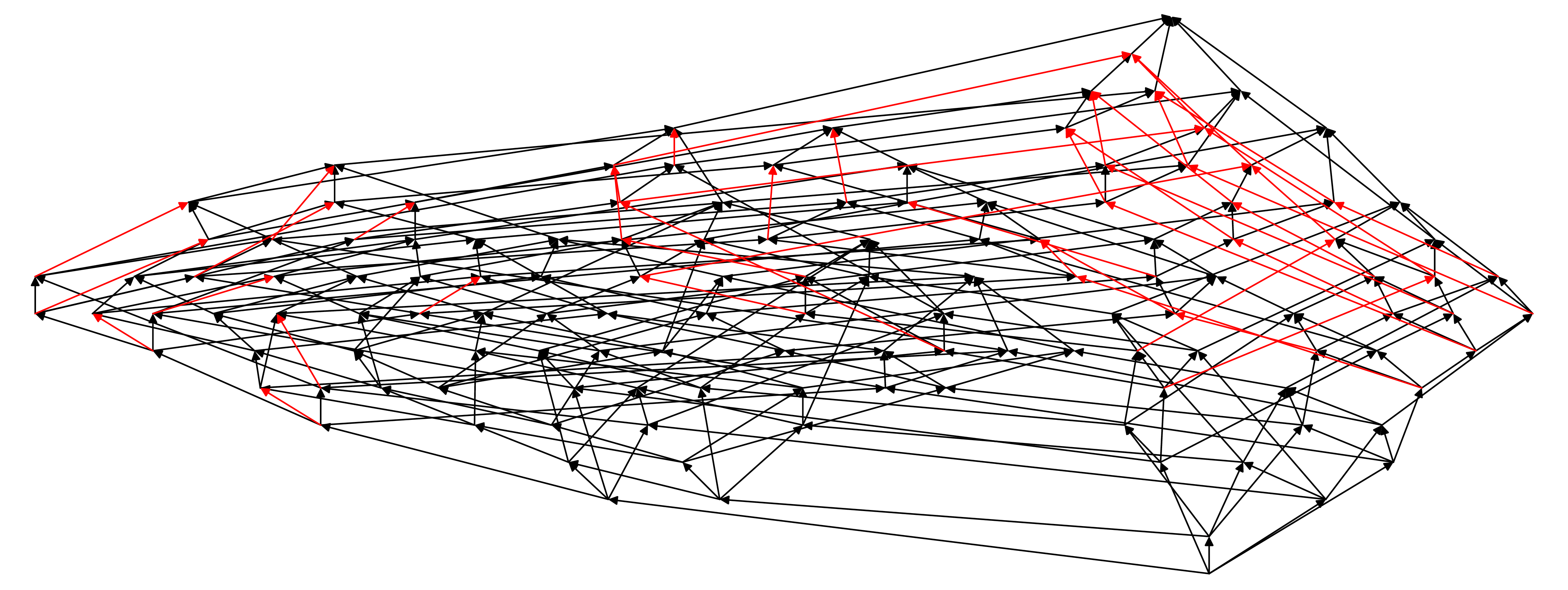}
     \caption{Hasse diagram of $\mathcal{S}_6=(\mathcal{D}_6,\leq)$ (black edges). The Tamari lattice $\mathcal{T}_6=(\mathcal{D}_6,\preceq)$ can be viewed by considering  black and red edges. }
     \label{fig6}
 \end{figure}

% https://www.irif.fr/~chapuy/fichiers/slides/LyonGdrIm2012.pdf

\section*{Acknowledgements} The authors are grateful to the referees for their detailed comments and corrections, which helped improve the paper. Jean-Luc Baril and Sergey Kirgizov were partially supported by ANR PiCs (ANR-22-CE48-0002) and ANR ARTICO funded by Bourgogne-Franche-Comt\'e region. Mehdi Naima was supported from the CNRS through the MITI interdisciplinary programs.

\end{document}